\documentclass[a4paper,11pt,twoside]{paper}
\usepackage[french,english]{babel}
\usepackage[latin1]{inputenc}
\usepackage{amsmath,amssymb,amsthm,mathrsfs,enumitem,tikz}
\usepackage{epic,eepic,pstricks}
\usepackage{hyperref}
\newtheoremstyle{sans}{\parskip}{\parskip}{\itshape}
                       {0pt}{\bfseries\sffamily}{.}{ }{}
\newtheoremstyle{sansplain}{\parskip}{\parskip}{}
                       {0pt}{\bfseries\sffamily}{.}{ }{}

\theoremstyle{sans}
\newtheorem{prop}{Proposition}[section]

\newtheorem{thm}[prop]{Theorem}
\newtheorem{lem}[prop]{Lemma}

\theoremstyle{sansplain}

\makeatletter
\@addtoreset{equation}{section}
\makeatother

\newcommand\1{\leavevmode\hbox{\rm \small1\kern-0.35em\normalsize1}}
      
\renewcommand{\geq}{\geqslant}
\renewcommand{\leq}{\leqslant}

\usepackage{color}
\definecolor{darkgreen}{rgb}{0,0.4,0}
\definecolor{MyDarkBlue}{rgb}{0,0.08,0.50}
\definecolor{BrickRed}{rgb}{0.65,0.08,0}

\hypersetup{
colorlinks=true,       
    linkcolor=blue,          
    citecolor=red,        
    filecolor=BrickRed,      
    urlcolor=darkgreen        
}

\begin{document}

\title{About a possible analytic approach for walks in the quarter plane with arbitrary big jumps\footnote{Version of \today}}

\author{Guy Fayolle\thanks{INRIA Paris-Rocquencourt, Domaine de Voluceau, BP 105, 78153 Le Chesnay Cedex, France. Email: \url{Guy.Fayolle@inria.fr}}    \and
        Kilian Raschel\thanks{CNRS \& F\'ed\'eration Denis Poisson \& Laboratoire de Math\'ematiques et Physique Th\'eorique, Universit\'e de Tours, Parc de Grandmont, 37200 Tours, France.
       Email: \url{Kilian.Raschel@lmpt.univ-tours.fr}}}

\maketitle

\begin{abstract}
In this note, we consider random walks in the quarter plane with arbitrary big jumps. We announce the extension to that class of models of the analytic approach of \cite{FIM}, initially valid for walks with small steps in the quarter plane. New  technical challenges arise, most of them being tackled in the framework of generalized boundary value problems on compact Riemann surfaces. 
\end{abstract}

\selectlanguage{francais}
\begin{abstract}
Dans cette note nous nous int\'eressons aux marches al\'eatoires avec sauts arbitrairement grands dans le quart de plan. Nous annon\c cons le d\'eveloppement pour cette classe de mod\`eles de l'approche analytique propos\'ee dans \cite{FIM}, initialement applicable aux marches \`a petits sauts dans le quart de plan. De nouvelles difficult\'es th\'eoriques surgissent, qui, pour l'essentiel, sont  abord\'ees dans le cadre de la th\'eorie des probl\`emes aux limites g\'en\'eralis\'es sur des surfaces de Riemann compactes. 
\end{abstract}
\selectlanguage{english}


\section{Introduction}
\label{}

In the past decades, many fruitful research activities have been dealing with the analysis of random walks in the quarter plane  (RWQP) or quadrant $\mathbb Z_+^2$. Indeed, these objects are at the crossroads of several domains. In our framework, the initial motivations were twofold:  on one hand, to analyze the stationary distribution of irreducible RWQP \cite{MAL}; on other hand, to study queueing models representing two coupled processors working at different service rates \cite{FI}. One can also consult \cite{CB,FIM} for modern reference books on analytic and probabilistic aspects of RWQP. Recently, applications were found in enumerative combinatorics, see \cite{BMM}. Indeed, walks in the quarter plane naturally encode many combinatorial objects (certain trees, maps, permutations, Young tableaux, etc.). They  also have many links with population biology and finance. Lastly, in another context, RWQP can be viewed as particular instances of random processes in cones. This latter topic is the subject of many recent works, due to its links with representation theory, quantum random walks, random matrices, non-colliding random walks, etc.

In most of the analytic studies (in particular in \cite{BMM,FI,FIM,MAL}), walks are supposed to have \emph{small steps}. This means that each unit of time, jumps take place to (a subset of) the eight nearest neighbors on the lattice, see Figure \ref{fig:Walks}. There are many reasons to make this hypothesis, their fundamental common point being that they  simplify the technical computations, and allow to get closed-form solutions. \emph{In this paper, we aim at presenting some challenges, and we announce results concerning the analysis of random walks having  possibly \emph{big jumps} in the quarter plane}. This extension of the analytic approach of \cite{FIM} has many possible applications, including queueing models, enumeration of lattice paths, discrete harmonic functions, stationary probabilities, etc.
\unitlength=1cm
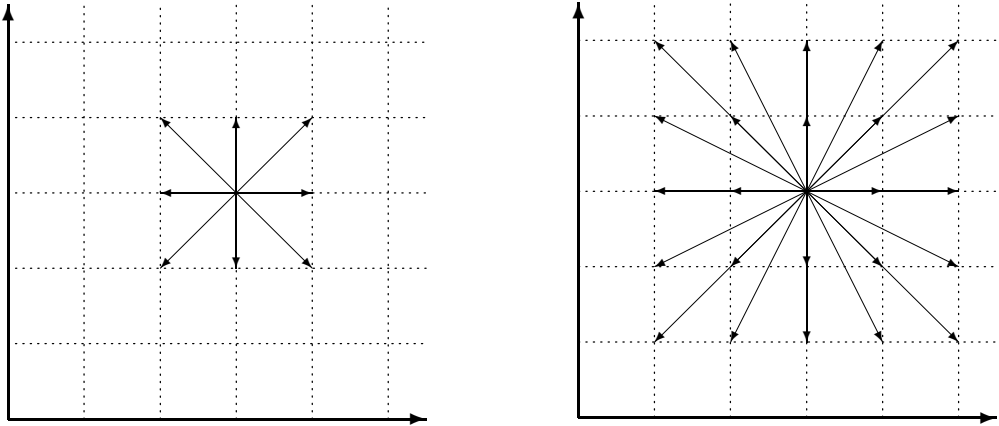
\begin{figure}[t]
\begin{center}
        \begin{picture}(17,7)
    \thicklines
    \thicklines
    \put(1,1){{\vector(1,0){5.5}}}
    \put(1,1){\vector(0,1){5.5}}
    \thinlines
    \put(4,4){\vector(1,1){1}}
    \put(4,4){\vector(-1,-1){1}}
    \put(4,4){\vector(1,0){1}}
    \put(4,4){\vector(-1,0){1}}
    \put(4,4){\vector(0,1){1}}
    \put(4,4){\vector(0,-1){1}}
    \put(4,4){\vector(-1,1){1}}
    \put(4,4){\vector(1,-1){1}}
    \linethickness{0.1mm}
    \put(1,2){\dottedline{0.1}(0,0)(5.5,0)}
    \put(1,3){\dottedline{0.1}(0,0)(5.5,0)}
    \put(1,4){\dottedline{0.1}(0,0)(5.5,0)}
    \put(1,5){\dottedline{0.1}(0,0)(5.5,0)}
    \put(1,6){\dottedline{0.1}(0,0)(5.5,0)}
    \put(2,1){\dottedline{0.1}(0,0)(0,5.5)}
    \put(3,1){\dottedline{0.1}(0,0)(0,5.5)}
    \put(4,1){\dottedline{0.1}(0,0)(0,5.5)}
    \put(5,1){\dottedline{0.1}(0,0)(0,5.5)}
    \put(6,1){\dottedline{0.1}(0,0)(0,5.5)}
\end{picture}

\vspace{-5mm}
        \begin{picture}(0,-10)
   \thicklines
    \put(1,1){{\vector(1,0){5.5}}}
    \put(1,1){\vector(0,1){5.5}}
    \thinlines
    \put(4,4){\vector(1,1){1}}
    \put(4,4){\vector(-1,-1){1}}
    \put(4,4){\vector(-1,-1){2}}
    \put(4,4){\vector(1,1){2}}
    \put(4,4){\vector(-1,0){2}}
    \put(4,4){\vector(1,0){2}}
    \put(4,4){\vector(0,1){2}}
    \put(4,4){\vector(0,-1){2}}
    \put(4,4){\vector(-1,1){2}}
    \put(4,4){\vector(1,-1){2}}
    \put(4,4){\vector(1,-2){1}}
    \put(4,4){\vector(1,2){1}}
    \put(4,4){\vector(-1,-2){1}}
    \put(4,4){\vector(-1,2){1}}
        \put(4,4){\vector(-2,1){2}}
    \put(4,4){\vector(2,1){2}}
    \put(4,4){\vector(-2,-1){2}}
    \put(4,4){\vector(2,-1){2}}
    \put(4,4){\vector(1,0){1}}
    \put(4,4){\vector(-1,0){1}}
    \put(4,4){\vector(0,1){1}}
    \put(4,4){\vector(0,-1){1}}
    \put(4,4){\vector(-1,1){1}}
    \put(4,4){\vector(1,-1){1}}
    \linethickness{0.1mm}
    \put(1,2){\dottedline{0.1}(0,0)(5.5,0)}
    \put(1,3){\dottedline{0.1}(0,0)(5.5,0)}
    \put(1,4){\dottedline{0.1}(0,0)(5.5,0)}
    \put(1,5){\dottedline{0.1}(0,0)(5.5,0)}
    \put(1,6){\dottedline{0.1}(0,0)(5.5,0)}
    \put(2,1){\dottedline{0.1}(0,0)(0,5.5)}
    \put(3,1){\dottedline{0.1}(0,0)(0,5.5)}
    \put(4,1){\dottedline{0.1}(0,0)(0,5.5)}
    \put(5,1){\dottedline{0.1}(0,0)(0,5.5)}
    \put(6,1){\dottedline{0.1}(0,0)(0,5.5)}
\end{picture}
\end{center}
  \vspace{-5mm}
\caption{On the left: walks with small jumps in the quarter plane; on the right: walks with arbitrary big jumps in the quarter plane (of maximal length $2$ on the figure)}
\label{fig:Walks}
\end{figure} 

To understand the intrinsic difficulties of such models, we should first recall that the now standard approach (for small steps walks) can be summarized by the three following steps.
\begin{enumerate}
\item\label{step:functional_equation} Find a functional equation for the generating function(s) of interest.
\item\label{step:BVP} Rewrite the functional equation as a \emph{boundary value problem} (BVP).
\item\label{step:solve_BVP} Solve the BVP.
\end{enumerate}

This approach concerns several types of problems pertaining to various mathematical areas, and allows to obtain many quantities: the stationary distribution of ergodic RWQP reflected on the boundary \cite{FIM,MAL}, the Green functions of killed random walks, the enumeration of deterministic walks \cite{BMM}, fine characteristics of certain queueing systems \cite{FI}, discrete harmonic functions, etc.

 Contrary to point (\ref{step:functional_equation}), we shall see that points (\ref{step:BVP}) and (\ref{step:solve_BVP}) raise more difficult issues, when assumptions on small jumps are relaxed. Indeed, point (\ref{step:functional_equation}) is rather simple, in the sense that most of the time, finding a functional equation simply reflects properties of the model. 
Point (\ref{step:BVP}), which first appeared in \cite{FI}, is the keystone of the whole approach.  
Point (\ref{step:solve_BVP}) is highly technical, and uses the standard literature (with some peculiarities) devoted to the resolution of BVPs.  In many examples, it turns out that the BVPs at stake (\ref{step:BVP}) have a unique solution, corresponding to the generating functions of interest. 


\section{Presentation of the model, functional equation, and principles of the approach}

\unitlength=1cm
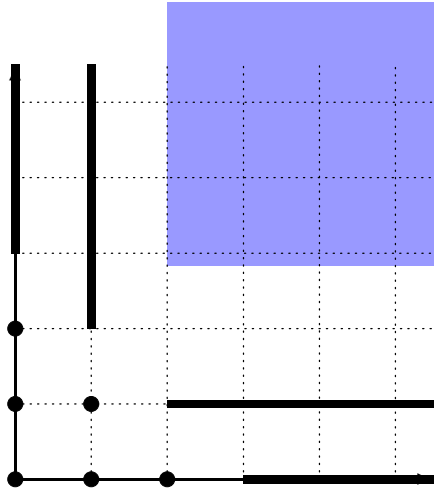
\begin{figure}[t]
\vspace{6mm}
\begin{center}
  \phantom{f}\hspace{20mm}\begin{tikzpicture}(19.55,6.5)
           \linethickness{2mm}
 \draw[draw=none,fill=blue!40] (10,3.5) -- (13.5,3.5) -- (13.5,0) -- (10,0) -- cycle;
     \end{tikzpicture}
     \end{center}
     \vspace{-35.5mm}
\begin{center}
     \begin{picture}(7.5,6.5)
    \thicklines
     \put(1,1){{\line(1,0){5.5}}}
    \put(1,1){\line(0,1){5.5}}
    \linethickness{1mm}
    \put(4,1){{\vector(1,0){2.5}}}
    \put(1,4){\vector(0,1){2.5}}
    \thinlines
\linethickness{1mm}
    \put(2,3){\line(0,2){3.5}}
    \put(3,2){\line(2,0){3.5}}
    \linethickness{0.1mm}
    \put(1,2){\dottedline{0.1}(0,0)(5.5,0)}
    \put(1,3){\dottedline{0.1}(0,0)(5.5,0)}
    \put(1,4){\dottedline{0.1}(0,0)(5.5,0)}
    \put(1,5){\dottedline{0.1}(0,0)(5.5,0)}
    \put(1,6){\dottedline{0.1}(0,0)(5.5,0)}
    \put(2,1){\dottedline{0.1}(0,0)(0,5.5)}
    \put(3,1){\dottedline{0.1}(0,0)(0,5.5)}
    \put(4,1){\dottedline{0.1}(0,0)(0,5.5)}
    \put(5,1){\dottedline{0.1}(0,0)(0,5.5)}
    \put(6,1){\dottedline{0.1}(0,0)(0,5.5)}
   {\put(2,2){\textcolor{black}{\circle*{0.2}}}}
   {\put(1,2){\textcolor{black}{\circle*{0.2}}}}
   {\put(1,3){\textcolor{black}{\circle*{0.2}}}}
   {\put(2,1){\textcolor{black}{\circle*{0.2}}}}
   {\put(3,1){\textcolor{black}{\circle*{0.2}}}}
   {\put(1,1){\textcolor{black}{\circle*{0.2}}}}
\end{picture}
\end{center}
  \vspace{-10mm}
\caption{If $I^-=J^-=2$ and $I^0=J^0=3$, there are $11$ domains of spatial homogeneity: the isolated points $(0,0)$, $(1,0)$, $(2,0)$, $(0,1)$, $(0,2)$ and $(1,1)$, the horizontal axes $\{(i,0) : i\geq 3\}$ and $\{(i,1) : i\geq 2\}$, the vertical axes $\{(0,j) : j\geq 3\}$ and $\{(1,j) : j\geq 2\}$, and the interior domain $\{(i,j) : i,j\geq 2\}$}
\label{fig:homogeneity}
\end{figure}

\subsection{The model}

In this section we introduce the notation, in the particular example of ergodic RWQP, although the announced results can be rendered much more general. 

Consider a two-dimensional random process with the following properties.
\begin{enumerate}[label={\rm (P\arabic{*})},ref={\rm (P\arabic{*})}]
\item\label{hypothesis:state_space}The state space is the quarter plane $\mathbb Z_+^2=\{0,1,2,\ldots \}^2$.
\item\label{hypothesis:space_homogeneity}The state space can be represented as the union of non-intersecting classes 
\begin{equation}
\label{eq:decomposition}
     \mathbb Z_+^2=S\cup\{\cup_\ell S'_{\ell}\}\cup \{\cup_k S''_{k}\}\cup\{\cup_{k,\ell} \{(k,\ell)\}\}.
\end{equation}
In order to describe the decomposition \eqref{eq:decomposition} we need to introduce four parameters, $(I^-,J^-)$ (resp.\ $(I^0,J^0)$), which describe the maximal negative amplitude of the transition probabilities in the interior (resp.\ on the boundary) of the quarter plane. Then, the interior class of the quadrant is
\begin{equation*}
     S=\{(i,j)\in\mathbb Z_+^2: i\geq I^-,j\geq J^-\}.
\end{equation*}
 In addition,
\begin{equation*} 
     S'_\ell=\{(i,\ell)\in\mathbb Z_+^2:i\geq I^-\},\qquad S''_k=\{(k,j)\in\mathbb Z_+^2:j\geq J^-\}
\end{equation*} 
are the horizontal and vertical domains strictly above the coordinate axes (i.e., $\ell\in\{1,\ldots ,J^--1\}$, $k\in\{1,\ldots ,I^--1\}$), and
\begin{equation*} 
     S'_0=\{(i,0)\in\mathbb Z_+^2:i\geq I^0\},\qquad S''_0=\{(0,j)\in\mathbb Z_+^2:j\geq J^0\}
\end{equation*} 
are parts of the coordinate axes. Lastly, we also need $\{(k,\ell)\}$ to take into account the finite number of remaining isolated points. See Figure \ref{fig:homogeneity}, where $I^-=J^-=2$ and $I^0=J^0=3$. Thus, in each class, the walk has state homogeneous transition probabilities, denoted by $p_{i,j}$, $p^{S'_\ell}_{i,j}$, $p^{S''_k}_{i,j}$ and $p^{(k,\ell)}_{i,j}$, respectively. 
\item\label{hypothesis:boundedness} We assume that $p_{i,j}=0  \mbox{ if } \{i<-I^- \mbox{ or } j<-J^-\} \mbox{ and if }\{i>I^+ \mbox{ or }j>J^+\}$. In this context, the small steps hypothesis reads $I^-=I^+=J^-=J^+=1$.
\end{enumerate}

\subsection{Statement of the functional equation}

We denote the generating function of the transition probabilities for the class $S$ by (when nothing is specified, the sum runs on all couples of indices $(i,j)\in\mathbb Z^2$)
\begin{equation}
\label{eq:def_R}
     R(x,y) = -1 + \textstyle \sum_{i,j}p_{i,j}x^iy^j ,
\end{equation}
and the stationary probabilities generating function by
\begin{equation}
\label{eq:def_pi}
     \pi(x,y) = \textstyle \sum_{i\geq I^-,j\geq J^-}\pi_{i,j}x^iy^j.
\end{equation}
For the other classes, we take the same notation with the associated class symbol in exponent (for instance, $R^{S'_\ell}(x,y)$ and $\pi^{S'_\ell}(x,y)$ for the class ${S'_\ell}$). Then the functional equation to be solved takes the form  (see \cite[Equation (1.3.4)]{FIM})
\begin{equation}
\label{eq:main_functional_equation}
     -R(x,y) \pi(x,y)  = \sum_{\{(k,\ell)\}}R^{(k,\ell)}(x,y)\pi^{(k,\ell)}(x,y)   
     +\sum_{\ell=0}^{J^--1}R^{S'_\ell}(x,y)\pi^{S'_\ell}(x,y)+\sum_{k=0}^{I^--1}R^{S''_k}(x,y)\pi^{S''_k}(x,y).
\end{equation}
It is valid (at least) in the bi-disc $\{x\in \mathbb C: \vert x\vert \leq 1\}\times\{x\in \mathbb C: \vert y\vert \leq 1\}$. The generating function $\pi^{S'_\ell}(x,y)$ (resp.\ $\pi^{S''_k}(x,y)$) is essentially a function of the single variable $x$ (resp.\ $y$), since $\pi^{S'_\ell}(x,y)=y^\ell \sum_{i\geq I^-} \pi_{i,\ell}x^i$ (resp.\ $\pi^{S''_k}(x,y)=x^k \sum_{j\geq J^-} \pi_{k,j}y^j$). Further, one has $\pi^{(k,\ell)}(x,y)=\pi_{k,\ell}x^ky^\ell$.

\subsection{Principles of the method}
\label{subsec:principles}
In the right-hand side of \eqref{eq:main_functional_equation}, there are $I^-+J^-$ unknown functions, together with a finite number of unknown coefficients $\pi_{k,\ell}$. We shall focus on the $I^-+J^-$ unknown functions, which all are sought a priori to be analytic  in the unit disc $\mathcal{D}$ (in fact, they will be only meromorphic in $\mathcal{D}$ when the system is not ergodic). The basic idea to determine these functions is to find $I^-+J^-$ independent equations. In our case these equations will take the form of BVPs set on closed curves (see \eqref{eq:boundary_condition} for an example), each of these problems being obtained from \emph{a pair of branch points in 
$\mathcal{D}$} of the polynomial 
\begin{equation}
\label{eq:def_P}
     P(x,y) = x^{I^-}y^{J^-}R(x,y),
\end{equation}
which is of degree $I=I^-+I^+$ in $x$ and $J=J^-+J^+$ in $y$, recalling that $x$ is a \emph{branch point} if the equation $P(x,y)=0$ in $y$ has a multiple root. In this probabilistic context, the unit disc condition is exactly the domain of definition (a priori) of the generating functions. From a technical point of view, a large part of the analysis will be devoted to locate these branch points, and to select which of them lead to boundary equations on closed curves, in connection with the underlying Riemann surface generated by the algebraic curve 
\begin{equation}
\label{eq:def_surface}
     \{(x,y)\in(\mathbb C\cup\{\infty\})^2 : P(x,y)=0\}.
\end{equation}


\section{Analytic approach for walks with small steps: branch points and the reduction to BVP}

Here we are left with only two unknowns in the right-hand side of the functional equation \eqref{eq:main_functional_equation}, as well as (to simplify) one unknown scalar coefficient $\pi^{(k,\ell)}(x,y)=\pi_{0,0}$ (in some sense, $\pi_{0,0}$ is not unknown, thanks to the normalization $\sum_{i,j\geq 0}\pi_{i,j}=1$). In the sequel, we shall write $\pi^{S'_\ell}(x,y)=y\pi(x)$ and $\pi^{S''_k}(x,y)=x\widetilde{\pi}(y)$.

\unitlength=0.6cm
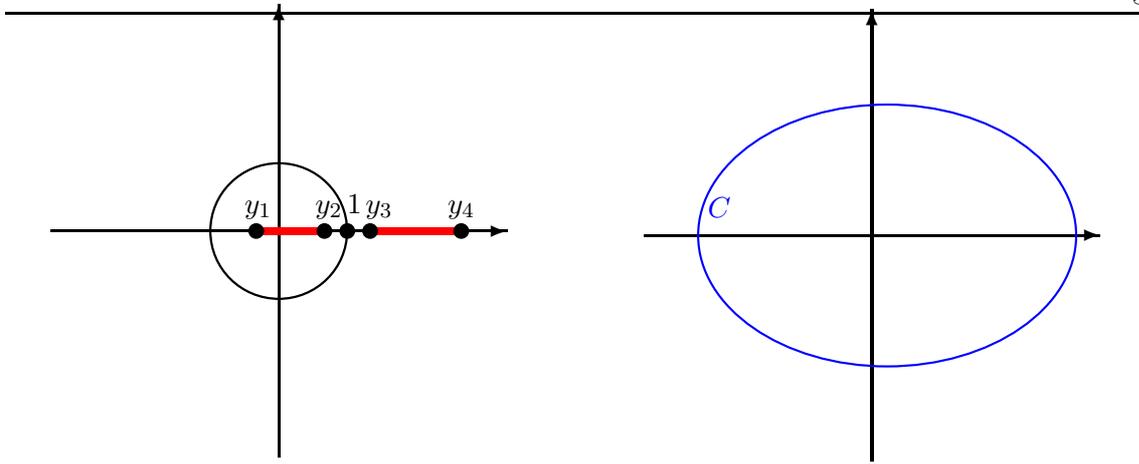
\begin{figure}[t]
\vspace{-40mm}
\begin{center}
    \begin{picture}(13,10)
    \thicklines
    \put(-5,0){\vector(1,0){10}}
    \put(0,-5){\vector(0,1){10}}
    \put(0,0){\circle{3}}
   \put(1.5,0.4){$1$}
   \put(-0.75,0.4){$y_1$}
   \put(0.8,0.4){$y_2$}
   \put(1.9,0.4){$y_3$}
   \put(3.7,0.4){$y_4$}
    {\put(1.5,0){\textcolor{black}{\circle*{0.3}}}}
    
\linethickness{1mm}
{\put(-0.5,0){\textcolor{red}{\line(1,0){1.5}}}}
{\put(2,0){\textcolor{red}{\line(1,0){2}}}}
{\put(-0.5,0){\textcolor{black}{\circle*{0.3}}}}
{\put(1,0){\textcolor{black}{\circle*{0.3}}}}
{\put(2,0){\textcolor{black}{\circle*{0.3}}}}
{\put(4,0){\textcolor{black}{\circle*{0.3}}}}
     \end{picture}
     
     \vspace{-4.2mm}
     
    \begin{picture}(-13,0)
    \thicklines
    \put(-5,0){\vector(1,0){10}}
    \put(0,-5){\vector(0,1){10}}
    \put(-3.6,0.4){\textcolor{blue}{$C$}}
    \psellipse[linecolor=blue](0.2,0)(2.5,1.75)
     \end{picture}
\end{center}
\vspace{27mm}
\caption{On the left: in the small steps case, there are two possible segments $L$ joining a couple of branch points, $[y_1,y_2]$ and $[y_3,y_4]$, but only one of them is included in the unit disc; on the right: a typical curve $C$ for the boundary condition \eqref{eq:boundary_condition}}
\label{fig:ex_curve}
\end{figure}

It was proved in \cite{FIM} that the functions $\pi(x)$ and $\widetilde \pi(y)$ can be found by solving a BVP on a closed curve. In more concrete terms, this means that there exists a curve $C$, symmetrical with respect to the horizontal axis (see Figure \ref{fig:ex_curve}), satisfying the boundary condition (below, we note $\alpha({t})=\overline{t}$ the complex conjugation, which is a particular automorphism of $C$, see Section \ref{sec:BC})
\begin{equation}
\label{eq:boundary_condition}
     A(t)\pi (t) - A(\alpha({t}))\pi(\alpha({t})) = g(t),\qquad \forall t\in C. 
\end{equation}
Above, $A(t)$ and $g(t)$ are known functions having simple expressions in terms of the parameters, see \cite[Equations (5.1.5) and (5.1.6)]{FIM}. The curve $C$ is obtained via a study of the branch points of the algebraic function of $y$ defined by $P(x,y)=0$. In the small steps case, with obvious notations, we may write $P(x,y)=b_2(y)x^2+b_1(y)x+b_0(y)$; its discriminant equals
\begin{equation*}
     \Delta(y)=b_1(y)^2-4b_2(y)b_0(y).
\end{equation*}
It is easily proved that this discriminant has  always $4$ real roots $y_k$ (in general not necessarily distinct), named branch points, as shown on Figure \ref{fig:ex_curve}. Let $L$ be a segment between two consecutive branch points. For $l\in L$, the two solutions to $P(x,l)=0$ are either real or complex conjugate. Choosing $L$ such that the two roots are complex conjugate, the set
\begin{equation}
\label{eq:expression_curve}
     \bigcup_{l\in L}\{x\in \mathbb C\cup\{\infty\} : P(x,l)=0\}
\end{equation}
is symmetrical with respect to the real axis in the $x$-plane, and it builds the two closed components of  a \emph{decomposed quartic curve}, which do not intersect and correspond respectively to the cuts $[y_1,y_2]$ and $[y_3,y_4]$ in the $x$-plane. On Figure \ref{fig:ex_curve}, $C$ is the component corresponding to $L=[y_1,y_2]$ in \eqref{eq:expression_curve}.

\section{Formal analytic approach for walks with arbitrary big steps}

We present a sequence of results showing to which extent the method of \cite{FIM} can be generalized to models of walks with arbitrary big jumps. The line of argument proceeds in stages, one of the key points being, as announced in Section \ref{subsec:principles}, to prove the existence of $I^-$ (resp.\ $J^-$) \emph{suitable} couples of branch points $x_k$ (resp.\ $y_k$) within the closed unit disc $\mathcal{D}$. Then each such pair will be associated with \emph{appropriate branches (or sheets)} of the algebraic curve \eqref{eq:def_surface} (after analytic continuation), and finally used to set a BVP for the $I^-$ (resp.\ $J^-$) unknown functions $\pi^{S''_k}(x,y)$ (resp.\ $\pi^{S'_\ell}(x,y)$). Needless to mention that selecting these so-called suitable branches involves some deep technicalities related to Riemann surfaces.    

\subsection{Finding and classifying the branch points}
Branch points in the $x$-plane are obtained by means of the \emph{discriminant} in $y$ of the polynomial
\begin{equation*}
     P(x,y)=b_I(y)x^I+\cdots+b_0(y)=a_J(x)y^J+\cdots+a_0(x).
 \end{equation*}

\begin{lem}\label{BP}
There are at most $2J(I-1)$ branch points $y_k$ and $2I(J-1)$ branch points $x_k$. 
\end{lem}
\smallskip
The proof of  Lemma \ref{BP} relies on the expression of the  resultant of $P(x,y)$ and $\frac{\partial P}{\partial x}(x,y)$ (resp.\ $\frac{\partial P}{\partial y}(x,y)$). These branch points can be split into two categories: the ones outside  $\mathcal{D}$, and those inside. We shall prove that, among the interior points, certain ones are not relevant for the problem, in the sense that they will not necessarily lead to boundary conditions on closed curves. This is the subject of  the next assertion.\smallskip

\begin{thm}
There are exactly $2I^-$ (resp.\ $2J^-$) interior branch points $x_k$ (resp.\ $y_k$)  leading to boundary conditions on closed curves. 
\end{thm}
\smallskip
The problem of leading or not to conditions on closed curves of the type \eqref{eq:boundary_condition}  strongly depends on the behavior of the branches of the algebraic functions defined by $P(x,Y(x))=P(X(y),y)=0$ in the neighborhood of the branch points. The answer to this difficult question is related to the genus and other specific properties of the underlying Riemann surface \eqref{eq:def_surface}.

\subsection{Branch cuts and boundary value problems}\label{sec:BC}
The $2J^-$ branch points $y_k$ can be grouped by pairs, forming segments $L$ in the $y$-plane. Some segments $L$ are real (as in the small steps case), but others can be complex conjugate (this is a new phenomenon pertaining to the big jumps case). For real segments, the boundary conditions are of  type \eqref{eq:boundary_condition}. But, for the non-real segments $L$,  condition \eqref{eq:boundary_condition}  changes: $\alpha$ is no more the complex conjugation, but a more general automorphism, sometimes called Carleman automorphism, satisfying $\alpha^2(t)=t$. \smallskip

\begin{thm}
The solution of the fundamental equation \eqref{eq:main_functional_equation} can be obtained from a generalized  BVP satisfied by a vector of $J^-$ unknown functions of the variable $x$. This BVP, defined  from  $J^-$ boundary conditions on closed curves, has in general a unique solution in the ergodic case. Mutatis mutandis, a similar statement holds for the $I^-$ unknown functions of $y$.  
\end{thm}
\smallskip
 In general, it is hopeless to expect a closed-form expression for the unknown functions in \eqref{eq:main_functional_equation}.  The reader will easily understand this claim, just by seeing the complexity of the computations in the small steps case \cite{FIM}. There is, however, a particular and interesting class for which we will be able to obtain explicit expressions, namely the models where $I^-=J^-=1$. Indeed, for these models, the functional equation \eqref{eq:main_functional_equation} takes the same form as for the small steps case. 

\subsection{Uniformization and {\L}ukasiewicz walks}
When the projections of the jumps in the west and south directions are of size one, while the others are arbitrary and bounded, one speaks sometimes of  \emph{{\L}ukasiewicz walks}, see \cite[Section I.5.3]{FS}. Then the following result holds.\smallskip
\begin{thm}
For  {\L}ukasiewicz walks, the unique unknown function $\pi(x)$ in the variable $x$ is given by 
\begin{equation*}
     \int_C\frac{f(t)w'(t)}{w(t)-w(x)}\textnormal{d}t.
\end{equation*}
Here, $f$ is a known algebraic function, and $w$ is directly expressed in terms of the conformal mapping between the interior domain of $C$ and the unit disc $\mathcal{D}$. This mapping in turn can be obtained by means of an ad hoc uniformization of the Riemann surface \eqref{eq:def_surface}.
\end{thm}
\medskip
The analysis of $w$ in the small steps case, done in \cite{FIM}, allows to introduce elliptic functions, since the Riemann surface \eqref{eq:def_surface} in the generic situation has genus one (the torus). But in the case of big jumps, the genus is bigger than one. Although the genus is not simply given in terms of the parameters $I^\pm$ and $J^\pm$, we can express it thanks to the Riemann-Hurwitz's formula, in terms of the order of the branch points.

\paragraph*{Acknowledgements}
The authors would like to thank the referee for his remarks about integrable systems and  possible applications to the quantum three-body problem.




\end{document}